\documentclass[12pt]{amsart}

\usepackage[all]{xy}
\usepackage{mathrsfs}
\usepackage{amsmath}
\usepackage{latexsym}
\usepackage{amssymb}
\usepackage{amsfonts}
\usepackage{longtable}
\usepackage{colortbl}
\include{PDF}
\usepackage{graphics}

\renewcommand{\proof}{\par\noindent{\it Proof.\ \ }}
\def\qed{\ifmmode\square\else\nolinebreak\hfill
$\Box$\fi\par\vskip12pt}

   \newcommand\Aut{\mathrm{Aut}}
 \newcommand\bbZ{\mathbb{Z}} 
\newcommand\C{\mathrm{C}} \newcommand\calB{\mathcal{B}} \newcommand\Cay{\mathrm{Cay}}

\newcommand\K{\mathsf{K}}

       \newcommand\Sy{\mathrm{S}} \newcommand\Sym{\mathrm{Sym}} 

 \newcommand\ZZ{{\mathbb Z}}

\def\bbZ{{\mathbb Z}}
\def\ov{\overline}

\def\s{\sigma}
\def\t{\tau}
\def\l{\langle}
\def\r{\rangle}

\def\Circ{{\bf Circ}}

\def\Ga{{\it \Gamma}}
\def\Sig{{\it \Sigma}}

\def\calB{{\mathcal B}}

\def\K{{\bf K}}

\newtheorem{theorem}{Theorem}[section]
\newtheorem{lemma}[theorem]{Lemma}

\newtheorem{corollary}[theorem]{Corollary}

\theoremstyle{definition}

\begin{document}

\title{Arc-transitive circulants}

\author[Song]{Shu Jiao Song}
\address{School of Mathematics and Science information, Yantai University}
\email{shujiao.song@ytu.edu.cn}
\thanks{2000 MR Subject Classification 20B15, 20B30, 05C25.}
\thanks{This work is partially supported by  NSFC (No. 61771019)}


\maketitle

\begin{abstract}
This short paper presents characterisations of normal arc-transitive circulants and arc-transitive normal circulants, that is, for a connected arc-transitive circulant $\Ga=\Cay(C,S)$, it is shown that
\begin{itemize}
\item $\Aut(C,S)$ is transitive on $S$ if and only if each element of $S$ has order $n$;
\item $\Aut\Ga\rhd C$ if and only if $S$ does not contain a coset of any subgroup.
\end{itemize}
This completes the classification of arc-transitive circulants given by Li-Xia-Zhou.
\end{abstract}

\section{Introduction}

A digraph is an ordered pair $(V,A)$ with the vertex set $V$ and arc set $A$,
in which an arc is an ordered pair of adjacent vertices.
The number $|V|$ of vertices is called the {\it order} of $\Gamma$.
A digraph $\Gamma$ is said to be \emph{arc-transitive} if the automorphism group $\Aut\Ga$ acts transitively on the arc set.

Let $C=\ZZ_n=\{0,1,\dots,n-1\}$ be a cyclic group of order $n$, and let $S$ be a subset of $C\setminus\{0\}$.
Let $\Cay(C,S)$ be the Cayley digraph with vertex set $C$ such that a vertex $v$ is adjacent to a vertex $w$ if and only if $w-v\in S$.
In this paper, we sometimes denote $\Cay(C,S)$ by $\Circ(n,S)$ if we do not need to emphasise on the connecting subset $S$.
A digraph is called a {\it circulant} of order $n$ if it is isomorphic to a digraph $\Circ(n,S)$ for some subset $S$ of $\ZZ_n$.

The cyclic group $C$ acts regularly on the vertex set by right multiplication as an automorphism group.
Thus a circulant is a vertex transitive graph.
A circulant $\Gamma$ is called a {\it normal circulant} if $\Aut\Ga$ has a cyclic subgroup
which is normal and regular on the vertices.

Much effort has been made to characterize arc-transitive circulants in the literature:
Chao and Wells~\cite{Chao1971,CW1973} classified those of prime order;
$2$-arc-transitive~\cite{ACMX1996,MW2000},
square-free order~\cite{LMM2001},
odd prime-power order~\cite{XBS2004}.
In the general case, a recursive characterization of arc-transitive circulants was given
by Kov\'{a}cs~\cite{Kovacs2004} and Li~\cite{Li2005}, independently.
Recently, Li-Xia-Zhou \cite{XiaBZ} prove that an arc-transitive circulant $\Ga$ has the form
\[\Ga\cong(\Ga_0\times\K_{n_1}\times\dots\times\K_{n_r})[\overline\K_b],\]
where $\Ga_0$ is a normal circulant of order $n_0$ and $n_0,n_1,\dots,n_r$ are pairwise coprime.
This would provide a complete characterization of arc-transitive circulants if normal circulants $\Ga_0$ is well-characterised.

For a Cayley graph $\Ga=\Cay(C,S)$, there is a useful subgroup of automorphisms 
\[\Aut(C,S)=\l\s\in\Aut(C)\mid S^\s=S\r.\]
This was first introduced by Godsil \cite{Godsil1981}.
Obviously, $\Aut\Ga\geqslant C{:}\Aut(C,S)$.
This motivated Xu \cite{Xu1998} to introduce {\it normal Cayley graphs}, which are those satisfying $\Aut\Ga=C{:}\Aut(C,S)$.
This was extended by Praeger \cite{Praeger1998} to define {\it normal arc-transitive Cayley graphs}, which are those such that $\Aut(C,S)$ is transitive on $S$.
We notice that a normal arc-transitive Cayley graph is not necessarily a normal Cayley graph.
Praeger \cite{Praeger1998} proposed to characterize normal arc-transitive Cayley graphs.

Obviously, if $\Aut(C,S)$ is transitive on $S$, then elements of $S$ all have the same order.
The first main result of this paper shows that the converse statement is also true, which is a bit surprising.

\begin{theorem}\label{thm-normal-arc-trans}
Let $\Ga=\Circ(n,S)$ be a connected arc-transitive circulant.
Then $\Ga$ is normal-arc-transitive if and only if each element of $S$ is of order $n$.
\end{theorem}

This solves a problem of \cite{Praeger1998} for cyclic groups.
Another version of Theorem~\ref{thm-normal-arc-trans} would be interesting.

\begin{corollary}
Let $\Ga=\Circ(n,S)$ be such that each element of $S$ has order $n$.
Then $\Ga$ is arc-transitive if and only if $\Ga$ is normal-arc-transitive.
\end{corollary}

In the terminology of \cite{XiaBZ}, normal arc-transitive circulants are characterized below, where $\pi(n)$ for an integer $n$ denotes the set of all prime divisors of $n$.

\begin{corollary}\label{normal-at}
A normal arc-transitive circulant has the form
\[\Ga=(\Ga_0\times\K_{p_1}\times\dots\times\K_{p_r})[\ov\K_b],\] 
such that $p_1,\dots,p_r$ are distinct primes, and $\pi(b)\subset\pi(n_0)\cup\{p_1,\dots, p_r\}$.
\end{corollary}

By definition, if an arc-transitive circulant is a normal circulant, then it is normal arc-transitive.
The next theorem provides a simple criterion for an arc-transitive circulant to be a normal circulant.

\begin{theorem}\label{thm-normal}
Let $\Ga=\Circ(n,S)$ be a connected arc-transitive circulant.
Then $\Ga$ is a normal circulant if and only if $S$ does not contain a coset of any subgroup.
\end{theorem}

\begin{corollary}\label{p-power}
An arc-transitive circulant $\Circ(n,S)$ of valency coprime to $n$ is a normal circulant.
In particular, $\Circ(p^e,S)$ with $p^e>3$ is normal if and only if $|S|$ divides $p-1$.
\end{corollary}

\begin{corollary}\label{(n,phi(n))-1}
Let $n=p_1^{e_1}\dots p_r^{e_r}$ be the prime factorisation.
If $\Circ(n,S)$ is a normal circulant, then the valency $|S|$ divides $(p_1-1)(p_2-1)\dots(p_r-1)$.
\end{corollary}

%
%

\section{Normal arc-transitive circulants}

Let $n$ be a positive integer, and let $\phi(n)$  the Euler $\phi$-function, which is the number of positive integers that are less than and coprime to $n$.
Let $C$ be a cyclic group of order $n$.
Then the automorphism group $\Aut(C)$ is an abelian group of order $\phi(n)$.


\begin{lemma}\label{prop3}
\emph{(\cite{Godsil1981,Xu1998})}
For a circulant $\Ga=\Cay(C,S)$,  the normalizer of $C$ in $\Aut\Ga$ is $C\rtimes\Aut(C,S)$.
\end{lemma}

A circulant $\Gamma$ is called {\it normal-arc-transitive}
if the normaliser of a cyclic regular subgroup in $\Aut\Ga$ is arc-transitive on $\Gamma$.

\begin{lemma}\label{quotient}
If $\Ga=\Sig[\ov\K_b]$ is a normal-edge-transitive circulant, then so is the quotient $\Sig$.
\end{lemma}

We note that not every arc-transitive circulant is normal-arc-transitive, for example, $\K_{p^2}$;
not every normal-arc-transitive circulant is a normal circulant,
for instance, $\K_p\times \K_q$ where $p,q$ are distinct primes.

\begin{lemma}\label{normal-arc-trans}
A circulant $\Circ(n,S)$ is normal arc-transitive if and only if $\l S\r=\ZZ_n$ and $\Aut(\ZZ_n,S)$ is regular on $S$.
\end{lemma}
\proof
Assume that $\Circ(n,S)$ is normal arc-transitive.
Then $\Aut(\ZZ_n,S)$ is transitive on $S$, and $\l S\r=\ZZ_n$.
Thus every element of $S$ is a generator of $\ZZ_n$.
If $\s\in\Aut(\ZZ_n,S)$ fixes some element $g\in S$, then since $\l g\r=\ZZ_n$, $\s$ fixes every element of $\ZZ_n$, and so $\s=$.
Therefore, $\Aut(\ZZ_n,S)$ is regular on $S$.

Conversely, if $\l S\r=\ZZ_n$ and $\Aut(\ZZ_n,S)$ is regular on $S$, then clearly $\Circ(n,S)$ is a normal arc-transitive circulant.
\qed

In other words, normal arc-transitive circulants are exactly normal arc-regular circulants.


For two graphs $\Ga=(U,E)$ and $\Sig=(V,F)$, the {\it direct product} $\Ga\times\Sig$ is the graph with vertex set $U\times V$ such that two vertices $(u_1,v_1),(u_2,v_2)$ are adjacent if and only if either $u_1,u_2$ are adjacent in $\Ga$ and $v_1=v_2$, or $v_1,v_2$ are adjacent in $\Sig$ and $u_1=u_2$.

\begin{lemma}\label{direct-prod}
Let $\Ga=\Cay(G_1,S_1)\times\Cay(G_2,S_2)$.
Then $\Ga=\Cay(G_1\times G_2,T)$, where $T=\{(s_1,s_2)\mid s_i\in S_i\}$, and furthermore, the following hold.
\begin{itemize}
\item[(i)] $\Ga$ is a circulant if and only if $G_1,G_2$ are cyclic groups of coprime orders.
\item[(ii)] Elements of $T$ have the same order if and only if elements of $S_i$ are all of the same order, where $i=1$ or $2$.
\end{itemize}
\end{lemma}
\proof
The proof follows from the definitions of Cayley graphs and direct product of graphs.
\qed

Obviously, a Cayley graph $\Cay(G,S)$ is a complete graph if and only if $S=G\setminus\{1\}$.

For digraphs $\Ga=(V_1,A_1)$ and $\Sig=(V_2,A_2)$, the \emph{lexicographic product} $\Ga[\Sig]$ of $\Sig$ by $\Ga$  is the digraph with vertex set $V_1\times V_2$ such that $(u_1,u_2)$ is connected to $(v_1,v_2)$ if and only if either $(u_1,v_1)\in A_1$, or $u_1=v_1$ and $(u_2,v_2)\in A_2$.

\begin{lemma}\label{key-1}
Let $C=\ZZ_n$ be a cyclic group of order $n$, and let $\Ga=\Cay(C,S)$ be connected and arc-transitive.
Then the following statements are equivalent:
\begin{itemize}
\item[(i)] $\Ga=\Sig[\ov\K_b]$ for some graph $\Sig$ and some integer $b>1$;
\item[(ii)] $S=T\l g_0\r$, where $g_0\in C$ and $T\subset S$, such that $|S|=|T||g_0|$ with $|g_0|=b$.
\end{itemize}
\end{lemma}
\proof
Assume first that $\Ga=\Sig[\ov\K_b]$.
Then $\Sig=\Cay(\ov C,\ov S)$ is a Cayley graph of a cyclic group $\ov C$, where $\ov C=C/N$ with $|N|=b$, and $|S|=|\ov S|b$.
It follows that $S=Ns_1\cup Ns_2\cup\dots\cup Ns_r$ such that $\ov S=\{\ov s_1,\ov s_2,\dots,\ov s_r\}$.
Let $N=\l g_0\r$.
Then $S=N\{s_1,s_2,\dots,s_r\}=\l g_0\r T$, where $T=\{s_1,s_2,\dots,s_r\}$.

Now assume that $S=T\l g_0\r$, where $T\subset S$.
Then
\[C=\l g_0\r\cup g\l g_0\r \cup \dots\cup g^{b-1}\l g_0\r ,\]
where $|g_0|=b$.
Let $B_i=\l g_0\r g^i$, where $0\leqslant i\leqslant b-1$.
Suppose $g^i\in S$ for some integer $i$.
Then $g^i\l g_0\r\subset S$, and so the vertex $(1,1)$ is adjacent to all vertices in $B_i=g^i\l g_0\r$.
It follows that the induced subgraph $[B_0,B_i]=\K_{b,b}$, and so $\Ga=\Sig[\ov\K_b]$, where $\Sig$ is the quotient graph $\Ga$ induced by the block system $\{B_0,B_1,\dots,B_{b-1}\}$.
\qed

We now quote the main result of \cite{XiaBZ}.

\begin{theorem}\label{arc-trans} {\rm(\cite{XiaBZ})}
Let $\Ga$ be a connected arc-transitive circulant of order $n$.
Then
\begin{itemize}
\item[(1)] $\Ga\cong(\Ga_0\times\K_{n_1}\times\dots\times\K_{n_r})[\overline\K_b]$, and
\item[(2)] $\Aut(\Ga)\cong\Sym(b)\wr\left(\Aut(\Ga_0)\times\Sym(n_1)\times\dots\times\Sym(n_r)\right)$,
\end{itemize}
where $\Ga_0\ncong\C_4$ is a connected normal circulant of order $n_0$, $n_i\geqslant4$ for $i=1,\dots,r$, and $n=n_0n_1\cdots n_rb$, and $n_0,n_1,\dots,n_r$ are pairwise coprime.
\end{theorem}

\begin{lemma}\label{key-0}
Let $\Ga=\Circ(n,S)$ be connected and arc-transitive.
If $S$ contains a coset of some subgroup, then $\Ga$ is not a normal circulant.
\end{lemma}
\proof
Suppose that $\Ga=\Cay(C,S)$ is a normal circulant, and $t\l h\r\subset S$ for some non-identity element $h\in C$
Then $\Aut(C,S)$ is transitive on $S$, and hence every element of $S$ is contained in a coset of $\l h\r$.
It follows that $S=t_1\l h\r\cup\dots\cup t_r\l h\r$.
By Lemma~\ref{key-1}, the graph $\Ga=\Sig[\ov\K_b]$, and $\Aut\Ga=\Sym(b)\wr\Aut\Sig$ by Theorem~\ref{arc-trans}.
So $\widehat C$ is not normal in $\Aut\Ga$.
\qed

\begin{lemma}\label{order-n}
Let $\Ga=\Circ(n,S)$ be such that all elements of $S$ have order $n$.
Then  $\Ga=(\Ga_0\times\K_{n_1}\times\dots\times\K_{n_r})[\ov\K_b]$ satisfying Theorem~$\ref{arc-trans}$ such that
\begin{itemize}
\item[(i)] $n_1,\dots,n_r$ are all primes,

\item[(ii)] $\pi(b)\subset\pi(n_0)\cup\{n_1,\dots, n_r\}$.
\end{itemize}
\end{lemma}
\proof
By Lemma~\ref{key-1}, we have $S=T\l h\r$ where $|h|=b$.
Let $T=\{t_1,\dots,t_r\}$ be such that $S$ is a disjoint union of cosets of $\l h\r$ with representatives in $T$:
\[S=t_1\l h\r\cup\dots \cup t_r\l h\r.\]
In particular, $t_1,\dots,t_r\in S$ since $\l h\r$ contains the identity.

Since elements of $S$ have the same order $n$, the elements $t_1,\dots,t_r$ are of order equal to $n$.
Suppose that a prime divisor $p$ of $b$ does not divide $n/b$.
If $|t_i|$ is coprime to $p$, then $t_ih$ has order divisible by $p$, which is a contradiction since $|t_i|=|t_ih|$.
Thus $t_i=xy$ where $|x|$ is a power of $p$ and $\gcd(|x|,|y|)=1$.
Then $x\in \l h\r$, and so $y=t_ix^{-1}\in S$, which is not possible since $|t_i|=|t_ix^{-1}|$.
So each prime divisor of $b$ divides $n/b$, namely, $\pi(b)\subset\pi(n/b)=\pi(n_0)\cup\{n_1,\dots,n_r\}$.
\qed

Let $C=\ZZ_n$, and let $n=p_1^{d_1}\dots p_t^{d_t}$ be the prime factorization of $n$.
Let $C_i$ be the Sylow $p_i$-subgroup of $C$.
Then $C_i=\ZZ_{p_i^{d_i}}$, and
\[\Aut(C_i)=
\left\{
\begin{array}{ll}
\ZZ_{p_i^{d_i-1}(p_i-1)}=\ZZ_{p_i^{d_i-1}}\times\ZZ_{p_i-1}, & \mbox{if $p_i$ is odd},\\
\ZZ_2\times\ZZ_{2^{d_i-2}}, &\mbox{if $p_i=2$.}
\end{array}
\right.
\]
In particular, the Hall $p_i'$-subgroup of $\Aut(C_i)$ is isomorphic to $\ZZ_{p_i-1}$, that is, by $\Aut(C_i)_{p_i'}=\ZZ_{p_i-1}$.
Moreover, $\Aut(C_1)_{p_1'}\times\dots\times\Aut(C_t)_{p_t'}$ induces a faithful action on the factor group $C/N$.


\begin{lemma}\label{Phi(C)}
Let $N\leqslant\Phi(C)$, and let $\ov C=C/N$.
Then each automorphism of $\ov C$ is induced by an automorphism of $C$, and
$\Aut(C)$ is surjective to $\Aut(\ov C)$.
\end{lemma}
\proof
It is easy to see that $\Aut(C)$ is regular on the set of generators of $C=\l g\r$.
Since $N$ is a characteristic subgroup of $C$, the automorphism group $\Aut(C)$ fixes $N$.
Now $\Aut(\ov C)$ acts regularly on the set of generators of $C/N=\l\ov g\r=\ZZ_m$.
Since $N\leqslant \Phi(C)$, we have $\pi(m)=\pi(n)$.
Let $\ov g^k$ and $\ov g^\ell$ be two generators of $\ov C$, and let $\ov\s\in\Aut(\ov C)$ such that $(\ov g^k)^{\ov\s}=\ov g^\ell$.
Then $\gcd(k,m)=\gcd(\ell,m)=1$, and so $\gcd(k,n)=\gcd(\ell,n)=1$.
Thus $g^k$ and $g^\ell$ are generators of $C$.
Let $\s\in\Aut(C)$ be such that $(g^k)^\s=g^\ell$.
Then $\s$ induces an automorphism of $\ov C$ which sends $\ov g^k$ to $\ov g^\ell$.
Clearly, this is the automorphism $\ov\s$ defined above.
\qed

\begin{lemma}\label{Lem-1}
Let $S=g\l h\r$ such that all elements of $S$ have order $n$ and $h\in N$ where $N\leqslant\Phi(C)$.
Then $\Aut(C,S)=\l \t\r$ is transitive on $S$, where $\t=(\s_1,\s_2,\dots,\s_t)$ is defined below.
\end{lemma}
\proof
Let $n=p_1^{d_1}p_2^{d_2}\dots p_t^{d_t}$, and let
\[g=g_1g_2\dots g_t\]
where $|g_i|=p_i^{d_i}$.
Then $h=h_1h_2\dots h_t$ such that $h_i=g_i^{e_i}$, where $e_i$ is a positive integer, and
$|h_i|=p_i^{d_i-e_i}$.
Let $S_i=g_i\l h_i\r$, for $1\leqslant i\leqslant t$.
Then
\[S_i=g_i\l h_i\r=\{g_i^{1+jp_i^{e_i}}\mid 0\leqslant j< d_i-e_i\}=\{g_i,g^{1+p_i^{e_i}},\dots,g^{1+?p_i^{e_i}}\}.\]
Let $\s_i\in\Aut(C_{p_i})$ be such that
\[g_i^{\s_i}=g_i^{1+p_i^{e_i}}.\]
Then $g_i^{\l\s_i\r}=\{g_i^{(1+p_i^{e_i})^k}\mid 0\leqslant k<d_i-e_i\}$.
Suppose that $g_i^{(1+p_i)^{k_1}}=g_i^{(1+p_i)^{k_2}}$ with $k_1>k_2$.
Then $g_i^{(1+p_i)^{k_1}-(1+p_i)^{k_2}}=1$, and so
\[(1+p_i)^{k_2}((1+p_i)^{k_1-k_2}-1)=0\ (\mod p_i^{d_i}).\]
So $(1+p_i)^{k_1-k_2}\equiv 1\ (\mod p_i^{d_i})$.
However, as $k_1-k_2<d_i$, this is not possible.
Thus the cardinality of $\{g_i^{(1+p_i^{e_i})^k}\mid 0\leqslant k<d_i-e_i\}$ is equal to $p_i^{e_i}$, and $\Aut(\l g_i\r,S_i)=\l\s_i\r$ is transitive on $S_i$.
Since $S=\{(s_1,s_2,\dots,s_t)\mid s_i\in S_i\}$ and $(\s_1,\s_2,\dots,\s_t)\in\Aut(C)$, it follows that  $\Aut(C,S)$ is transitive on $S$.
\qed

\begin{lemma}\label{normal-cir}
An arc-transitive circulant $\Cay(C,S)$ is a normal circulant if and only if $\Aut(P_i)_{p_i}\cap\Aut(C,S)=1$ for $1\leqslant i\leqslant t$ and $\Aut(P_i)_{p_i'}\not\leqslant\Aut(C,S)$ if $P_i=\ZZ_{p_i}$.

\end{lemma}
\proof
Suppose that $\Aut(P_i)_{p_i}\cap\Aut(C,S)=\l s_1\r\not=1$.
Then for each $g\in S$, we have $g^{\l\s_1\r}=g\l h\r$ where $|h|=|\s_1|$.
So $S=T\l h\r$, and $\Ga=(\Ga_0\times\K_{n_1}\times\dots\times\K_{n_r})[\ov\K_b]$ with $b=|h|$, satisfying the conditions given in Theorem~\ref{arc-trans}.
Thus $\Ga$ is not a normal circulant by Theorem~\ref{arc-trans}.
If $P_i=\ZZ_{p_i}$ and $\Aut(P_i)_{p_i'}\not\leqslant\Aut(C,S)$, then $\Ga=\Sig\times\K_{p_i}$, which implies that $\Ga$ is not a normal circulant.

If $\Aut(P_i)_{p_i}\cap\Aut(C,S)=1$ for $1\leqslant i\leqslant t$ and $\Aut(P_i)_{p_i'}\not\leqslant\Aut(C,S)$ whenever $P_i=\ZZ_{p_i}$, then by Theorem~\ref{arc-trans} we have $\Ga=\Ga_0$, which is a normal circulant.
\qed

\begin{lemma}\label{key-2}
Let $\Ga=\Circ(C,S)$ be such that all elements of $S$ have order equal to $n$.
Then $\Aut(C,S)$ is transitive on $S$, and $\Ga$ is a normal arc-transitive circulant.
\end{lemma}
\proof
Since elements of $S$ are all of order $n$, by Lemma~\ref{order-n}, $\Ga=(\Ga_0\times\K_{n_1}\times\dots\times\K_{n_r})[\ov\K_b]$ satisfying Theorem~$\ref{arc-trans}$ is such that $n_1,\dots,n_r$ are all distinct primes, and $\pi(b)\subset\pi(n_0)\cup\{n_1,\dots, n_r\}$.
By Lemma~\ref{key-1}, the connecting set $S$ has the form $S=T\l h\r$, where $|h|=b$.
Let $p_i=n_i$ for $1\leqslant i\leqslant r$, and
let $\pi(n_0)=\{p_{r+1},\dots,p_{r+s}\}$ if $n_0\not=1$.

Let $\Ga_0=\Cay(n_0,T_0)$, and $\K_{p_i}=\Cay(p_i,T_i)$, where $T_i=\ZZ_{p_i}\setminus\{1\}$.
Then
\[T=\{(t_0,t_1,\dots,t_r)\mid t_i\in T_i\ \mbox{for}\ 0\leqslant i\leqslant r\}.\]
Thus
\[S=T\l h\r={\Large\cup}_{(t_0,t_1,\dots,t_r)\in T}(t_0,t_1,\dots,t_r)\l h\r.\]

By Theorem~\ref{arc-trans}, the automorphism group
\[\Aut(\Ga)\cong\Sym(b)\wr\left(\Aut(\Ga_0)\times\Sym(p_1)\times\dots\times\Sym(p_r)\right)\]
where $p_1,\dots,p_r$ are distinct primes.
For $1\leqslant i\leqslant r+s$, let $b_{p_i}$ be the $p_i$-part of $b$.
Then $\l h_{p_i}\r$ is a Sylow $p_i$-subgroup of $\l h\r$.
Since $p_i$ divides $n/b$, there exists an automorphism $\s_i$ of a Sylow $p_i$-subgroup $C_{p_i}$ of $C$ of order equal to $|h_{p_i}|$.
For $1\leqslant i\leqslant r$, let $\t_i$ be an automorphism of $C_{p_i}$ of order $p_i-1$.

Since $\Ga_0$ is a normal circulant, $\Ga_0=\Cay(C_0,S_0)$ with $C_0=\ZZ_{n_0}$ such that $\Aut(C_0,S_0)$ is transitive on $S_0$.
As $\pi(b)\subset\pi(n_0)\cup\{p_1,\dots, p_r\}$, each prime $p_i$ divides $n_0$ for $r+1\leqslant i\leqslant r+s$.
Let $C_{\pi(n_0)}$ be the Hall $\pi(n_0)$-subgroup of $C$.
Then by Lemma~\ref{Frattini} there exists a subgroup $K\leqslant\Aut(C_{p_i})$ such that $\Aut(C_0,S_0)$ is the quotient of $K$ and $K\cong\Aut(C_0,S_0)$.

Now $\Ga_\calB=\Ga_0\times\K_{p_1}\times\dots\times\K_{p_r}=\Cay(\ov C,\ov S)$, and $X/X_{(\calB)}$ is arc-transitive on $\Ga_\calB$.
Clearly, all elements of $\ov S$ have the same order, and by induction, $\Aut(\ov C,\ov S)$ is transitive on $\ov S=\{\ov g_1,\dots,\ov g_r\}$.
For any $j\in\{1,\dots,r\}$, there exists $\ov\s_j\in\Aut(\ov C,\ov S)$ such that
\[\ov g_1^{\ov\s_j}=\ov g_j.\]
By Lemma~\ref{Phi(C)}, there exists $\s_j\in\Aut(C)$ which is a preimage of $\ov\s_j$.
Then $(g_1\l h\r)^{\s_j}=g_j\l h\r$, and so
\[g_1^{\s_j}=g_jh^{\xi_j}\]
for some integer $\xi_j$.
Let $\t$ be the automorphism defined above, and let
\[H=\l \t,\s_2,\dots,\s_r\r.\]
Then $H\leqslant\Aut(C,S)$ and $H$ is transitive on $S$.
So $\Ga$ is a normal arc-transitive circulant.
\qed

\vskip0.1in
{\bf Proof of Theorem~\ref{thm-normal-arc-trans}:}

Assume that $\Ga=\Cay(C,S)$ is a normal edge-transitive circulant, where $C=\ZZ_n$.
Then $\Aut(C,S)$ is transitive on $S$, and hence all elements of $S$ have the same order.
On the other hand, by Theorem~\ref{arc-trans}, $\Ga\cong(\Ga_0\times\K_{n_1}\times\dots\times\K_{n_r})[\overline\K_b]$.
By Lemma~\ref{quotient}, the quotient $\Sig=\Cay(\ov C,\ov S)$ is a norma-edge-transitive circulant.
If $n_i$ is not a prime for $i\geqslant 1$, then $\ov S$ contains elements of different orders, not possible.
Thus each $n_i$ is a prime, as in part~(ii).

Suppose part (ii) holds, namely, $\Ga\cong(\Ga_0\times\K_{p_1}\times\dots\times\K_{p_r})[\overline\K_b]$, where
$n_0=|\Ga_0|$ and either $b=1$ or $\pi(b)\subset\pi(n_0p_1\dots p_r)$ such that  $p_1,\dots,p_r$ are distinct primes, $n=n_0p_1\dots p_r b$, and $\gcd(n_0,p_1\dots p_r)=1$.
Further, if $\Ga_0\not=\C_4$, then $\Aut(\Ga)\cong\Sy_b\wr(\Aut(\Ga_0)\times\Sy_{p_1}\times\dots\times\Sy_{p_r})$.

\section{Normal circulants}

Let $\Ga=\Cay(C,S)$ be connected and arc-transitive, where $C$ is cyclic.

Assume that $\Ga$ is a normal circulant, namely, $\Aut\Ga=C{:}\Aut(C,S)$.
Then $\Aut(C,S)$ is transitive on $S$.
By Lemma~\ref{key-0}, $S$ does not contain any coset $g\l h\r$.
Suppose that $S$ contains $g(\l h\r^\#)$.
Then 
\[S=\{g_1,\dots,g_r\}\l h\r^\#=g_1\l h\r^\#\cup\dots\cup g_r\l h\r^\#\]
since $\Aut(C,S)$ is transitive on $S$ and $\l h\r$ is a characteristic subgroup of $C$.

Let $C=\l g\r=\ZZ_n$ and $n=p_1^{e_1}\dots p_r^{e_r}$ be a prime decomposition.
Let $C_i$ be the Sylow $p_i$-subgroup of $C$.
Then $\Aut(C)=\Aut(C_1)\times\dots\times\Aut(C_r)$.

Assume that $p_1<p_2<\dots<p_t$, and $P_i$ is the Sylow $p_i$-subgroup of $C$.
Then $C=P_1\times P_2\times\dots\times P_t$ with $|P_i|=p_i^{e_i}$, and
\[\Aut(C)=\Aut(P_1)\times\Aut(P_2)\times\dots\times\Aut(P_t),\]
and $\Aut(P_i)=\ZZ_{(p_i-1)p_i^{e_i-1}}$ if $p_i$ is odd, and $\Aut(\ZZ_{2^e})=\ZZ_2\times\ZZ_{2^{e-2}}$ for $e\geqslant2$.

\begin{lemma}\label{normal-cir}
An arc-transitive circulant $\Cay(C,S)$ is a normal circulant if and only if $\Aut(C,S)\leqslant\ZZ_{p_1-1}\times\ZZ_{p_2-1}\times\dots\times\ZZ_{p_t-1}$.
\end{lemma}
\proof
Assume that $\Ga=\Cay(C,S)$ is a normal Cayley graph which is arc-transitive.
Then $\Cay(C,S)$ is an arc-transitive normal circulant, and by Theorem~\ref{Arc-trans-normal},
\[\Ga=(\Ga_0\times\K_{p_1}\times\dots\times\K_{p_r})[\ov\K_b]\]
with $n_0=|\Ga_0|$ such that either $b=1$ or $\pi(b)\subset\pi(n_0p_1\dots p_r)$,
where $p_1,\dots,p_r$ are distinct primes, $n=n_0p_1\dots p_r b$, and $\gcd(n_0,p_1\dots p_r)=1$.

{\bf Proof of Corollary~\ref{p-power}:}\ \
Let $\Circ(p^e,S)$ be a normal circulant and arc-transitive.
Then it is neither a complete graph nor a lexicographic product of a circulant and an isolated graph.
Thus $\gcd(|S|,p)=1$.
Since the automorphism group of a cyclic group $\ZZ_{p^e}$ is a cyclic group of order $p^e(p-1)$,
it follows that $|S|$ divides $p-1$.

Conversely, assume that $\Circ(p^e,S)$ is arc-transitive and $|S|$ divides $p-1$ and less then $p^e-1$.
Let $\Ga=\Circ(p^e,S)$, and let $G=\Aut\Ga$.

Since $|S|$ divides $p-1$, we have $\Ga\not=\Sig[\ov\K_b]$.
Since $|S|<p^e-1$, $\Ga$ is not a complete graph.
Since $|G|=p^e$, it follows that $\Ga\not=\Ga_0\times\dots\times\Ga_m$, where $m\geqslant1$.
Thus $\Circ(p^e,S)=\Ga_0$ is a normal circulant.
\qed

Let $\Ga=\Circ(n,S)$.
If $b\not=1$, then
\[S=S'B,\]
where $S'$ is subset of $S$ and $B$ is a subgroup of $R=\ZZ_n$.
If $b=1$ and $\Ga=\Ga_1\times\K_p$ with $p$ prime, then
\[S=T\ZZ_p^\#,\]
where $T$ is the projection of $S$ in $R_{p'}$.

\begin{lemma}\label{direct-decomp}
A circulant is unitary if and only if it is a direct product of its unitary Sylow-circulants.
\end{lemma}
\proof
Let $\Gamma=\Circ(n,S)$ be a unitary circulant of order $n$.
Then $\Gamma$ has a direct product decomposition:
\[\Circ(n,S)=\Circ({p_1^{e_1}},S_1)\times\dots\times\Circ({p_r^{e_r}},S_r),\]
where each $\Circ(p_i^{e_i},S_i)$ is a unitary circulant.
\qed

A circulant is normal arc-transitive if and only if all of its Sylow projections are normal arc-transitive.

The unit circulant $\Circ(n,S)$ is such that $S$ consists of all generators
of $\bbZ_n$.

\begin{corollary}
For $n=p_1^{e_1}\dots p_r^{e_r}$, the unit circulant of order $n$ is isomorphic to
\[(\K_{p_1}\times\dots\times\K_{p_r})[\ov\K_{p_1^{e_1-1}\dots p_r^{e_r-1}}]=
\K_{p_1}[\ov\K_{p_1^{e_1-1}}]\times\dots\times\K_{p_r}[\ov\K_{p_r^{e_r-1}}].\]
\end{corollary}

\vskip0.1in

\end{document}